\documentclass[a4paper,12pt]{amsart}

\usepackage{booktabs}
\usepackage{amsmath, amsfonts, amssymb, amsthm}
\usepackage[all]{xy}
\usepackage[height=22.5cm, width=15.5cm]{geometry}

\numberwithin{equation}{section}

\theoremstyle{plain}
\newtheorem{thm}{Theorem}[section]
\newtheorem{prop}[thm]{Proposition}
\newtheorem{lem}[thm]{Lemma}
\newtheorem{cor}[thm]{Corollary}
\newtheorem*{thm*}{Theorem}

\theoremstyle{definition}

\theoremstyle{remark}
\newtheorem*{rem}{Remark}

\newcommand{\mbb}[1]{\mathbb{#1}}
\newcommand{\wt}[1]{\widetilde{#1}}

\newcommand{\ol}[1]{\overline{#1}}

\newcommand{\abs}[1]{\lvert #1\rvert}

\newcommand{\hq}{/\hspace{-0.12cm}/}

\DeclareMathOperator{\im}{Im}

\DeclareMathOperator{\Aut}{Aut}

\DeclareMathOperator{\eps}{\varepsilon}

\DeclareMathOperator{\Sing}{Sing}

\providecommand{\bysame}{\leavevmode ---\ }
\providecommand{\og}{``}
\providecommand{\fg}{''}
\providecommand{\smfandname}{\&}

\title{On proper $\mbb{R}$-actions on hyperbolic Stein surfaces}

\author{Christian Miebach}
\address{LATP-UMR(CNRS) 6632, CMI-Universit\'e d'Aix-Marseille I, 39, rue
Joliot-Curie, F-13453 Marseille Cedex 13, France. }
\email{miebach@cmi.univ-mrs.fr}

\author{Karl Oeljeklaus}
\address{LATP-UMR(CNRS) 6632, CMI-Universit\'e d'Aix-Marseille I, 39, rue
Joliot-Curie, F-13453 Marseille Cedex 13, France.}
\email{karloelj@cmi.univ-mrs.fr}

\thanks{The authors would like to thank Peter Heinzner and Jean-Jacques Loeb
for numerous discussions on the subject.}

\begin{document}

\begin{abstract}
In this paper we investigate proper $\mbb{R}$--actions on hyperbolic Stein
surfaces and prove in particular the following result: Let $D\subset\mbb{C}^2$
be a simply-connected bounded domain of holomorphy which admits a proper
$\mbb{R}$--action by holomorphic transformations. The quotient $D/\mbb{Z}$ with
respect to the induced proper $\mbb{Z}$--action is a Stein manifold. A normal
form for the domain $D$ is deduced.
\end{abstract} 

\maketitle 

\section{Introduction}

Let $X$ be a Stein manifold endowed with a real Lie transformation group $G$ of
holomorphic automorphisms. In this situation it is natural to ask whether there
exists a $G$--invariant holomorphic map $\pi\colon X\to X\hq G$ onto a complex
space $X\hq G$ such that $\mathcal{O}_{X\hq G}=(\pi_*\mathcal{O}_X)^G$ and, if
yes, whether this quotient $X\hq G$ is again Stein. If the group $G$ is
compact, both questions have a positive answer as is shown in~\cite{He}.

For non-compact $G$ even the existence of a complex quotient in the above sense
of $X$ by $G$ cannot be guaranteed. In this paper we concentrate on the most
basic and already non-trivial case $G=\mbb{R}$. We suppose that $G$ acts
properly on $X$. Let $\Gamma=\mbb{Z}$. Then $X/\Gamma$ is a complex manifold
and if, moreover, it is Stein, we can define $X\hq G:= (X/\Gamma) \hq
(G/\Gamma)$. The following was conjectured by Alan~Huckleberry.

\begin{quote}
{\emph{Let $X$ be a contractible bounded domain of holomorphy in $\mbb{C}^n$
with a proper action of $G=\mbb{R}$. Then the complex manifold $X/\Gamma$ is
Stein.}}
\end{quote}

In~\cite{FI} this conjecture is proven for the unit ball and in~\cite{Mie4} for
arbitrary bounded homogeneous domains in $\mbb C^n$. In this paper we make a
first step towards a proof in the general case by showing

\begin{thm*}
Let $D$ be a simply-connected bounded domain of holomorphy in $\mbb{C}^2$.
Suppose that the group $\mbb{R}$ acts properly by holomorphic transformations
on $D$. Then the complex manifold $D/\mbb{Z}$ is Stein. Moreover, $D/\mbb{Z}$
is biholomorphically equivalent to a domain of holomorphy in $\mbb{C}^2$.
\end{thm*}

As an application of this theorem we deduce a normal form for domains of
holomorphy whose identity component of the automorphism group is non-compact as
well as for proper $\mbb{R}$--actions on them. Notice that we make no
assumption on smoothness of their boundaries.

We first discuss the following more general situation. Let $X$ be a hyperbolic
Stein manifold with a proper $\mbb{R}$--action. Then there is an induced local
holomorphic $\mbb{C}$--action on $X$ which can be globalized in the sense
of~\cite{HeIan}. The following result is central for the proof of the above
theorem.

\begin{thm*}
Let $X$ be a hyperbolic Stein surface with a proper $\mbb{R}$--action. Suppose
that either $X$ is taut or that it admits the Bergman metric and
$H^1(X,\mbb{R})=0$. Then the universal globalization $X^*$ of the induced local
$\mbb{C}$--action is Hausdorff and $\mbb{C}$ acts properly on $X^*$.
Furthermore, for simply-connected $X$ one has that $X^*\to X^*/\mbb{C}$ is a
holomorphically trivial $\mbb{C}$--principal bundle over a simply-connected
Riemann surface.
\end{thm*}

Finally, we discuss several examples of hyperbolic Stein manifolds $X$ with
proper $\mbb{R}$--actions such that $X/\mbb{Z}$ is not Stein. If one does not
require the existence of an $\mbb{R}$--action, there are bounded Reinhardt
domains in $\mbb{C}^2$ with proper $\mbb{Z}$--actions for which the quotients
are not Stein.

\section{Hyperbolic Stein $\mbb{R}$--manifolds}

In this section we present the general set-up.

\subsection{The induced local $\mbb{C}$--action and its globalization}

Let $X$ be a hyperbolic Stein manifold. It is known that the group $\Aut(X)$ of
holomorphic automorphisms of $X$ is a real Lie group with respect to the
compact-open topology which acts properly on $X$ (see~\cite{Kob2}). Let
$\{\varphi_t\}_{t\in\mbb{R}}$ be a closed one parameter subgroup of $\Aut(D)$.
Consequently, the action $\mbb{R}\times X\to X$, $t\cdot x:=\varphi_t(x)$, is
proper. By restriction, we obtain also a proper $\mbb{Z}$--action on $X$. Since
every such action must be free, the quotient $X/\mbb{Z}$ is a complex manifold.
This complex manifold $X/\mbb{Z}$ carries an action of $S^1\cong\mbb{R}/
\mbb{Z}$ which is induced by the $\mbb{R}$--action on $X$.

Integrating the holomorphic vector field on $X$ which corresponds to this
$\mbb{R}$--action we obtain a local $\mbb{C}$--action on $X$ in the following
sense. There are an open neighborhood $\Omega\subset\mbb{C}\times X$ of $\{0\}
\times X$ and a holomorphic map $\Phi\colon\Omega\to X$, $\Phi(t,x)=:t\cdot x$,
such that the following holds:
\begin{enumerate}
\item For every $x\in X$ the set $\Omega(x):=\bigl\{t\in\mbb{C};\ (t,x)\in
\Omega\bigr\}\subset\mbb{C}$ is connected;
\item for all $x\in X$ we have $0\cdot x=x$;
\item we have $(t+t')\cdot x=t\cdot(t'\cdot x)$ whenever both sides are
defined.
\end{enumerate}

Following~\cite{Pa} (compare \cite{HeIan} for the holomorphic setting) we say
that a globalization of the local $\mbb{C}$--action on $X$ is an open
$\mbb{R}$--equivariant holomorphic embedding $\iota\colon X\hookrightarrow X^*$
into a (not necessarily Hausdorff) complex manifold $X^*$ endowed with a
holomorphic $\mbb{C}$--action such that $\mbb{C}\cdot \iota(X)=X^*$. A
globalization $\iota\colon X\hookrightarrow X^*$ is called {\it universal} if
for every $\mbb{R}$--equivariant holomorphic map $f\colon X\to Y$ into a
holomorphic $\mbb{C}$--manifold $Y$ there exists a holomorphic
$\mbb{C}$--equivariant map $F\colon X^*\to Y$ such that the diagram
\begin{equation*}
\xymatrix{
X\ar[rr]^{\iota}\ar[dr]_f & & X^*\ar[dl]^F\\
& Y &
}
\end{equation*}
commutes. It follows that a universal globalization is unique up to
isomorphism if it exists.

Since $X$ is Stein, the universal globalization $X^*$ of the induced local
$\mbb{C}$--action exists as is proven in~\cite{HeIan}. We will always identify
$X$ with its image $\iota(X)\subset X^*$. Then the local $\mbb{C}$--action on
$X$ coincides with the restriction of the global $\mbb{C}$--action on $X^*$ to
$X$.

Recall that $X$ is said to be orbit-connected in $X^*$ if for every $x\in X^*$
the set $\Sigma(x):=\{t\in\mbb{C};\ t\cdot x\in X\}$ is connected. The
following criterion for a globalization to be universal is proven
in~\cite{IanTaTr}.

\begin{lem}\label{Lem:Orbit-conn}
Let $X^*$ be any globalization of the induced local $\mbb{C}$--action on $X$.
Then $X^*$ is universal if and only if $X$ is orbit-connected in $X^*$.
\end{lem}

\begin{rem}
The results about (universal) globalizations hold for a bigger class of groups
(\cite{IanTaTr}). However, we will need it only for the groups $\mbb{C}$ and
$\mbb{C}^*$ and thus will not give the most general formulation.
\end{rem}

For later use we also note the following

\begin{lem}
The $\mbb{C}$--action on $X^*$ is free.
\end{lem}

\begin{proof}
Suppose that there exists a point $x\in X^*$ such that $\mbb{C}_x$ is
non-trivial. Because of $\mbb{C}\cdot X=X^*$ we can assume that $x\in X$ holds.
Since $\mbb{C}_x$ is a non-trivial closed subgroup of $\mbb{C}$, it is either a
lattice of rank $1$ or $2$, or $\mbb{C}$. The last possibility means that $x$
is a fixed point under $\mbb{C}$ which is not possible since $\mbb{R}$ acts
freely on $X$.

We observe that the lattice $\mbb{C}_x$ is contained in the connected
$\mbb{R}$--invariant set $\Sigma(x)=\{t\in\mbb{C};\ t\cdot x\in X\}$. By
$\mbb{R}$--invariance $\Sigma(x)$ is a strip. Since $X$ is hyperbolic, this
strip cannot coincide with $\mbb{C}$. The only lattice in $\mbb{C}$ which can
possibly be contained in such a strip is of the form $\mbb{Z}r$ for some $r\in
\mbb{R}$. Since this contradicts the fact that $\mbb{R}$ acts freely on $X$,
the lemma is proven.
\end{proof}

Note that we do not know whether $X^*$ is Hausdorff. In order to guarantee the
Hausdorff property of $X^*$, we make further assumptions on $X$. The following
result is proven in~\cite{Ian2} and~\cite{IanSpTr}.

\begin{thm}\label{Thm:Existence}
Let $X$ be a hyperbolic Stein manifold with a proper $\mbb{R}$--action. Suppose
in addition that $X$ is taut or admits the Bergman metric. Then $X^*$ is
Hausdorff. If $X$ is simply-connected, then the same is true for $X^*$.
\end{thm}

We refer the reader to Chapter~4.10 and Chapter~5 in~\cite{Kob2} for the
definitions and examples of tautness and the Bergman metric.

\begin{rem}
Every bounded domain in $\mbb C^n$ admits the Bergman metric.
\end{rem}

\subsection{The quotient $X/\mbb{Z}$}

We assume from now on that $X$ fulfills the hypothesis of
Theorem~\ref{Thm:Existence}. Since $X^*$ is covered by the translates $t\cdot
X$ for $t\in\mbb{C}$ and since the action of $\mbb{Z}$ on each domain $t\cdot
X$ is proper, we conclude that the quotient $X^*/\mbb{Z}$ fulfills all axioms
of a complex manifold except for possibly not being Hausdorff.

We have the following commutative diagram:
\begin{equation*}
\xymatrix{
X\ar[r]\ar[d] & X^*\ar[d]\\
X/\mbb{Z}\ar[r] & X^*/\mbb{Z}.
}
\end{equation*}
Note that the group $\mbb{C}^*=(S^1)^\mbb{C}\cong\mbb{C}/\mbb{Z}$ acts on
$X^*/\mbb{Z}$. Concretely, if we identify $\mbb{C}/\mbb{Z}$ with $\mbb{C}^*$
via $\mbb{C}\to\mbb{C}^*$, $t\mapsto e^{2\pi it}$, the quotient map $p\colon
X^*\to X^*/\mbb{Z}$ fulfills $p(t\cdot x)=e^{2\pi it}\cdot p(x)$.

\begin{lem}
The induced map $X/\mbb{Z}\hookrightarrow X^*/\mbb{Z}$ is the universal
globalization of the local $\mbb{C}^*$--action on $X/\mbb{Z}$.
\end{lem}

\begin{proof}
The open embedding $X\hookrightarrow X^*$ induces an open embedding $X/\mbb{Z}
\hookrightarrow X^*/\mbb{Z}$. This embedding is $S^1$--equivariant and we have
$\mbb{C}^*\cdot X/\mbb{Z}=X^*/\mbb{Z}$. This implies that $X^*/\mbb{Z}$ is a
globalization of the local $\mbb{C}^*$--action on $X/\mbb{Z}$.

In order to prove that this globalization is universal,  
by the globalization theorem in \cite{IanTaTr} it is enough to show that
$X/\mbb{Z}$ is orbit-connected in $X^*/\mbb{Z}$. Hence, we must show that for
every $[x]\in X/\mbb{Z}$ the set $\Sigma\bigl([x]\bigr):=\{t\in\mbb{C}^*;\
t\cdot[x]\in X/\mbb{Z}\}$ is connected in $\mbb{C}^*$. For this we consider the
set $\Sigma(x)=\{t\in\mbb{C};\ t\cdot x\in X\}$. Since the map $X\to X/\mbb{Z}$
intertwines the local $\mbb{C}$-- and $\mbb{C}^*$--actions, we conclude that
$t\in\Sigma(x)$ holds if and only if $e^{2\pi it}\in\Sigma\bigl([x]\bigr)$
holds. Since $X^*$ is universal, $\Sigma(x)$ is connected which implies that
$\Sigma\bigl([x]\bigr)$ is likewise connected. Thus $X^*/\mbb{Z}$ is universal.
\end{proof}

\begin{rem}
The globalization $X^*/\mbb{Z}$ is Hausdorff if and only if $\mbb{Z}$ or,
equivalently, $\mbb{R}$ act properly on $X^*$. As we shall see in
Lemma~\ref{proper}, this is the case if $X$ is taut.
\end{rem}

\subsection{A sufficient condition for $X/\mbb{Z}$ to be Stein}

If $\dim X=2$, we have the following sufficient condition for $X/\mbb{Z}$ to be
a Stein surface.

\begin{prop}\label{Prop:suffcrit}
If the $\mbb{C}$--action on $X^*$ is proper and if the Riemann surface
$X^*/\mbb{C}$ is not compact, then $X/\mbb{Z}$ is Stein.
\end{prop}

\begin{proof}
Under the above hypothesis we have the $\mbb{C}$--principal bundle $X^*\to
X^*/\mbb{C}$. If the base $X^*/\mbb{C}$ is not compact, then this bundle is
holomorphically trivial, i.\,e.\ $X^*$ is biholomorphic to $\mbb{C}\times R$
where $R$ is a non-compact Riemann surface. Since $R$ is Stein, the same is
true for $X^*$ and for $X^*/\mbb{Z}\cong\mbb{C}^*\times R$. Since $X/\mbb{Z}$
is locally Stein, see \cite{Mie4}, in the Stein manifold $X^*/\mbb{Z}$, the
claim follows from~\cite{DocGr}.
\end{proof}

Therefore, the crucial step in the proof of our main result consists in showing
that $\mbb{C}$ acts properly on $X^*$ under the assumption $\dim X=2$.

\section{Local properness}

Let $X$ be a hyperbolic Stein $\mbb{R}$--manifold. Suppose that $X$ is taut or
that it admits the Bergman metric and $H^1(X,\mbb{R})=\{0\}$. We show that then
$\mbb{C}$ acts locally properly on $X^*$.

\subsection{Locally proper actions}

Recall that the action of a Lie group $G$ on a manifold $M$ is called
locally proper if every point in $M$ admits a $G$--invariant open neighborhood
on which the $G$ acts properly.

\begin{lem}
Let $G\times M\to M$ be locally proper.
\begin{enumerate}
\item For every $x\in M$ the isotropy group $G_x$ is compact.
\item Every $G$--orbit admits a geometric slice.
\item The orbit space $M/G$ is a smooth manifold which is in general not
Hausdorff.
\item All $G$--orbits are closed in $M$.
\item The $G$--action on $M$ is proper if and only if $M/G$ is Hausdorff.
\end{enumerate}
\end{lem}

\begin{proof}
The first claim is elementary to check. The second claim is proven
in~\cite{DuiKo}. The third one is a consequence of $(2)$ since the slices yield
charts on $M/G$ which are smoothly compatible because the transitions are given
by the smooth action of $G$ on $M$. Assertion $(4)$ follows from $(3)$ because
in locally Euclidian topological spaces points are closed. The last claim is 
proven in~\cite{Pa2}.
\end{proof}

\begin{rem}
Since $\mbb R$ acts properly on $X$, the $\mbb R$-action on $X^*$ is
locally proper.
\end{rem}

\subsection{Local properness of the $\mbb{C}$--action on $X^*$}

Recall that we assume that
\begin{equation}\label{Eqn:taut}
\text{$X$ is taut}
\end{equation}
or that
\begin{equation}\label{Eqn:Bergman}
\text{$X$ admits the Bergman metric and $H^1(X,\mbb{R})=\{0\}$.}
\end{equation}

We first show that assumption~\eqref{Eqn:taut} implies that $\mbb{C}$ acts
locally properly on $X^*$. 

Since $X^*$ is the universal globalization of the induced local
$\mbb{C}$--action on $X$, we know that $X$ is orbit-connected in $X^*$. This
means that for every $x\in X^*$ the set $\Sigma(x)=\{t\in\mbb{C};\ t\cdot x\in
X\}$ is a strip in $\mbb{C}$. In the following we will exploit the properties
of the thickness of this strip.

Since $\Sigma(x)$ is $\mbb{R}$--invariant, there are ``numbers''
$u(x)\in\mbb{R}\cup\{-\infty\}$ and $o(x)\in\mbb{R}\cup\{\infty\}$ for every
$x\in X^*$ such that
\begin{equation*}
\Sigma(x)=\bigl\{t\in\mbb{C};\ u(x)<\im(t)<o(x)\bigr\}.
\end{equation*}

The functions $u\colon X^*\to\mbb{R}\cup\{-\infty\}$ and
$o\colon X^*\to\mbb{R}\cup\{\infty\}$ so obtained are upper and lower
semicontinuous, respectively. Moreover, $u$ und $o$ are $\mbb{R}$--invariant
and $i\mbb{R}$--equivariant:
\begin{equation*}
u(it\cdot x)=u(x)-t\quad\text{and}\quad o(it\cdot x)=o(x)-t.
\end{equation*}

\begin{prop}
The functions $u,-o\colon X^*\to\mbb{R}\cup\{-\infty\}$ are plurisubharmonic.
Moreover, $u$ and $o$ are continuous on $X^*\setminus\{u=-\infty\}$ and
$X^*\setminus\{o=\infty\}$, respectively.
\end{prop}

\begin{proof}
It is proven in \cite{For} that $u$ and $-o$ are plurisubharmonic
on $X$. By equivariance, we obtain this result for $X^*$.

Now we prove that the function $u\colon X\setminus\{u=-\infty\}\to\mbb{R}$ is
continuous which was remarked without complete proof in \cite{Ian2}.
For this let $(x_n)$ be a sequence in $X$ which converges to
$x_0\in X\setminus\{u=-\infty\}$. Since $u$ is upper semi-continuous, we have
$\limsup_{n\to\infty}u(x_n)\leq u(x_0)$. Suppose that $u$ is not continuous in
$x_0$. Then, after replacing $(x_n)$ by a subsequence, we find $\eps>0$ such
that $u(x_n)\leq u(x_0)-\eps<u(x_0)$ holds for all $n\in\mbb{N}$. Consequently,
we have $\Sigma(x_0) =\bigl\{t\in\mbb{C};\ u(x_0)<\im(t)<o(x_0)\bigr\}
\subset\Sigma:=\bigl\{t\in \mbb{C};\ u(x_0)-\eps<\im(t)<o(x_0)\bigr\}
\subset\Sigma(x_n)$ for all $n$ and hence obtain the sequence of holomorphic
functions $f_n\colon\Sigma\to X$, $f_n(t):=t\cdot x_n$. Since $X$ is taut and
$f_n(0)=x_n\to x_0$, the sequence $(f_n)$ has a subsequence which compactly
converges to a holomorphic function $f_0\colon\Sigma\to X$. Because of
$f_0\bigl(iu(x_0)\bigr)=\lim_{n\to\infty}f_n\bigl(iu(x_0)\bigr)=\lim_{n\to
\infty}iu(x_0)\cdot x_n=iu(x_0)\cdot x_0\notin X$ we arrive at a contradiction.
Thus the function $u\colon X\setminus\{u=-\infty\}\to\mbb{R}$ is continuous. By
$(i\mbb{R})$--equivariance, $u$ is also continuous on
$X^*\setminus\{u=-\infty\}$. A similar argument shows continuity of $-o\colon
X^*\setminus\{o=\infty\}\to\mbb{R}$.
\end{proof}

Let us consider the sets
\begin{equation*}
{\mathcal{N}}(o):=\bigl\{x\in X^*;\ o(x)=0\bigr\}\quad\text{and}\quad
{\mathcal{P}}(o):=\bigl\{x\in X^*;\ o(x)=\infty\bigr\}.
\end{equation*}
The sets ${\mathcal{N}}(u)$ and ${\mathcal{P}}(u)$ are similarly defined. Since
$X=\bigl\{x\in X^*;\ u(x)<0<o(x)\bigr\}$, we can recover $X$ from $X^*$ with
the help of $u$ and $o$.

\begin{lem}\label{proper}
The action of $\mbb{R}$ on $X^*$ is proper.
\end{lem}

\begin{proof}
Let $\partial^*X$ denote the boundary of $X$ in $X^*$. Since the functions $u$
and $-o$ are continuous on $X^*\setminus\mathcal{P}(u)$ and
$X^*\setminus\mathcal{P}(o)$ one verifies directly that
$\partial^*X=\mathcal{N}(u)\cup\mathcal{N}(o)$ holds. As a consequence, we note
that if $x\in\partial^*X$, then for every $\eps>0$ the element $(i\eps)\cdot x$
is not contained in $\partial^*X$.

Let $(t_n)$ and $(x_n)$ be sequences in $\mbb{R}$ and $X^*$ such that
$(t_n\cdot x_n,x_n)$ converges to $(y_0,x_0)$ in $X^*\times X^*$. We may assume
without loss of generality that $x_0$ and hence $x_n$ are contained in $X$ for
all $n$. Consequently, we have $y_0\in X\cup\partial^*X$. If
$y_0\in\partial^*X$ holds, we may choose an $\eps>0$ such that $(i\eps)\cdot
y_0$ and $(i\eps)\cdot x_0$ lie in $X$. Since the $\mbb{R}$--action on $X$ is
proper, we find a convergent subsequence of $(t_n)$ which was to be shown.
\end{proof}

\begin{lem}
We have:
\begin{enumerate}
\item ${\mathcal{N}}(u)$ and ${\mathcal{N}}(o)$ are $\mbb{R}$--invariant.
\item We have ${\mathcal{N}}(u)\cap{\mathcal{N}}(o)=\emptyset$.
\item The sets ${\mathcal{P}}(u)$ and ${\mathcal{P}}(o)$ are closed,
$\mbb{C}$--invariant and
pluripolar in $X^*$.
\item ${\mathcal{P}}(u)\cap{\mathcal{P}}(o)=\emptyset$.
\end{enumerate}
\end{lem}

\begin{proof}
The first claim follows from the $\mbb{R}$--invariance of $u$ and $o$.

The second claim follows from $u(x)<o(x)$.

The third one is a consequence of the $\mbb{R}$--invariance and
$i\mbb{R}$--equivariance of $u$ and $o$.

If there was a point $x\in\mathcal{P}(u)\cap\mathcal{P}(o)$, then $\mbb{C}\cdot
x$ would be a subset of $X$ which is impossible since $X$ is hyperbolic.
\end{proof}

\begin{lem}\label{Lem:Slice}
If $o$ is not identically $\infty$, then the map
\begin{equation*}
\varphi \colon i\mbb{R}\times{\mathcal{N}}(o)\to
X^*\setminus{\mathcal{P}}(o),\quad \varphi(it,z)=
it\cdot z,
\end{equation*}
is an $i\mbb{R}$--equivariant homeomorphism. Since $\mbb{R}$ acts properly on
${\mathcal{N}}(o)$, it follows that $\mbb{C}$ acts properly on
$X^*\setminus{\mathcal{P}}(o)$. The same holds when $o$ is replaced by $u$.
\end{lem}

\begin{proof}
The inverse map $\varphi^{-1}$ is given by $x\mapsto\bigl(-io(x),io(x)\cdot
x\bigr)$.
\end{proof}

\begin{cor}
The $\mbb{C}$--action on $X^*$ is locally proper. If
${\mathcal{P}}(o)=\emptyset$ or ${\mathcal{P}}(u)=\emptyset$ hold, then
$\mbb{C}$ acts properly on $X^*$.
\end{cor}

From now on we suppose that $X$ fulfills the assumption~\eqref{Eqn:Bergman}.
Recall that the Bergman form $\omega$ is a K\"ahler form on $X$ invariant under
the action of $\Aut(X)$. Let $\xi$ denote the complete holomorphic vector field
on $X$ which corresponds to the $\mbb{R}$--action, i.\,e.\ we have
$\xi(x)=\left.\frac{\partial}{\partial t}\right|_0\varphi_t(x)$. Hence,
$\iota_\xi\omega=\omega(\cdot,\xi)$ is a $1$--form on $X$ and since
$H^1(X,\mbb{R})=\{0\}$ there exists a function
$\mu^\xi\in\mathcal{C}^\infty(X)$ with $d\mu^\xi=\iota_\xi\omega$.

\begin{rem}
This means that $\mu^\xi$ is a momentum map for the $\mbb{R}$--action on $X$.
\end{rem}

\begin{lem}
The map $\mu^\xi\colon X\to\mbb{R}$ is an $\mbb{R}$--invariant submersion.
\end{lem}

\begin{proof}
The claim follows from $d\mu^\xi(x)J\xi_x=\omega_x(J\xi_x,\xi_x)>0$.
\end{proof}

\begin{prop}
The $\mbb{C}$--action on $X^*$ is locally proper.
\end{prop}

\begin{proof}
Since $\mu^\xi$ is a submersion, the fibers $(\mu^\xi)^{-1}(c)$, $c\in\mbb{R}$,
are real hypersurfaces in $X$. Then
\begin{equation*}
\left.\frac{d}{dt}\right|_0\mu^\xi(it\cdot x)=\omega_x(J\xi_x,\xi_x)>0
\end{equation*}
implies that every $i\mbb{R}$--orbit intersects $(\mu^\xi)^{-1}(c)$
transversally. Since $X$ is orbit-connected in $X^*$, the map
$i\mbb{R}\times(\mu^\xi)^{-1}(c)\to X^*$ is injective and therefore a
diffeomorphism onto its open image. Together with the fact that
$(\mu^\xi)^{-1}(c)$ is $\mbb{R}$--invariant this yields the existence of
differentiable local slices for the $\mbb{C}$--action.
\end{proof}

\subsection{A necessary condition for $X/\mbb{Z}$ to be Stein}

We have the following necessary condition for $X/\mbb{Z}$ to be a Stein
manifold.

\begin{prop}\label{Prop:Proper}
If the quotient manifold $X/\mbb{Z}$ is Stein, then $X^*$ is Stein and
the $\mbb{C}$--action on $X^*$ is proper.
\end{prop}

\begin{proof}
Suppose that $X/\mbb{Z}$ is a Stein manifold. By~\cite{IanTaTr} this implies
that $X^*$ is Stein as well.

Next we will show that the $\mbb{C}^*$--action on $X^*/\mbb{Z}$ is proper. For
this we will use as above a moment map for the $S^1$--action on $X^*/\mbb{Z}$.

By compactness of $S^1$ we may apply the complexification theorem
from~\cite{He} which shows that $X^*/\mbb{Z}$ is also a Stein manifold and in
particular Hausdorff. Hence, there exists a smooth strictly plurisubharmonic
exhaustion function $\rho\colon X^*/\mbb{Z}\to\mbb{R}^{>0}$ invariant under
$S^1$. Consequently, $\omega:=\frac{i}{2}\partial\overline{\partial}
\rho\in\mathcal{A}^{1,1}(X^*)$ is an $S^1$--invariant K\"ahler form. Associated
to $\omega$ we have the $S^1$--invariant moment map
\begin{equation*}
\mu\colon X^*/\mbb{Z}\to\mbb{R},\quad\mu^\xi(x):=\left.\frac{d}{dt}\right|_0
\rho\bigl(\exp(it\xi)\cdot x\bigr),
\end{equation*}
where $\xi$ is the complete holomorphic vector field on $X^*/\mbb{Z}$ which
corresponds to the $S^1$--action. Now we can apply the same argument as above
in order to deduce that $\mbb{C}^*$ acts locally properly on $X^*/\mbb{Z}$.

We still must show that $(X^*/\mbb{Z})/\mbb{C}^*$ is Hausdorff. To see this,
let $\mbb{C}^*\cdot x_j$, $j=0,1$, be two different orbits in $X^*/\mbb{Z}$.
Since $\mbb{C}^*$ acts locally properly, these are closed and therefore there
exists a function $f\in\mathcal{O}(X^*/\mbb{Z})$ with $f|_{\mbb{C}^*\cdot
x_j}=j$ for $j=0,1$. Again we may assume that $f$ is $S^1$-- and consequently
$\mbb{C}^*$--invariant. Hence, there is a continuous function on
$(X^*/\mbb{Z})/\mbb{C}^*$ which separates the two orbits, which implies that
$(X^*/\mbb{Z})/\mbb{C}^*$ is Hausdorff. This proves that $\mbb{C}^*$ acts
properly on $X^*/\mbb{Z}$.

Since we know already that the $\mbb{C}$--action on $X^*$ is locally proper, it
is enough to show that $X^*/\mbb{C}$ is Hausdorff. But this follows from the
properness of the $\mbb{C}^*$--action on $X^*/\mbb{Z}$ since $X^*/\mbb{C}\cong
(X^*/\mbb{Z})/\mbb{C}^*$ is Hausdorff.
\end{proof}

\section{Properness of the $\mbb{C}$--action}

Let $X$ be a hyperbolic Stein $\mbb{R}$--manifold. Suppose that $X$
fulfills~\eqref{Eqn:taut} or~\eqref{Eqn:Bergman}. We have seen that $\mbb{C}$
acts locally properly on $X^*$. In this section we prove that under the
additional assumption $\dim X=2$ the orbit space $X^*/\mbb{C}$ is Hausdorff.
This implies that $\mbb{C}$ acts properly on $X^*$ if $\dim X=2$.

\subsection{Stein surfaces with $\mbb{C}$--actions}

For every function $f\in\mathcal{O}(\Delta)$ which vanishes only at the origin,
we define
\begin{equation*}
X_f:=\bigl\{(x,y,z)\in\Delta\times\mbb{C}^2;\ f(x)y-z^2=1\bigr\}.
\end{equation*}
Since the differential of the defining equation of $X_f$ is given by $\bigl(
f'(x)y\ f(x)\ -2z\bigr)$, we see that $1$ is a regular value of $(x,y,z)\mapsto
f(x)y-z^2$. Hence, $X_f$ is a smooth Stein surface in $\Delta\times\mbb{C}^2$.

There is a holomorphic $\mbb{C}$--action on $X_f$ defined by
\begin{equation*}
t\cdot(x,y,z):=\bigl(x,y+2tz+t^2f(x),z+tf(x)\bigr).
\end{equation*}
One can directly check that this defines an action.

\begin{lem}
The $\mbb{C}$--action on $X_f$ is free, and all orbits are closed.
\end{lem}

\begin{proof}
Let $t\in\mbb{C}$ such that $(x,y+2tz+t^2f(x),z+tf(x)\bigr)=(x,y,z)$ for
some $(x,y,z) \in X_f$. If $f(x)\not=0$, then $z+tf(x)=z$ implies $t=0$. If
$f(x)=0$, then $z\not=0$ and $y+2tz=y$ gives $t=0$.

The map $\pi\colon X_f\to\Delta$, $(x,y,z)\mapsto x$, is $\mbb{C}$--invariant.
If $a\in\Delta^*$, then $f(a)\not=0$ and we have
\begin{equation*}
\frac{z}{f(a)}\cdot\bigl(a,f(a)^{-1},0\bigr)=(a,y,z)\in X_f,
\end{equation*}
which implies $\pi^{-1}(a)=\mbb{C}\cdot\bigl(a, f(a)^{-1},0\bigr)$. A similar
calculation gives $\pi^{-1}(0)=\mbb{C}\cdot p_1\cup\mbb{C}\cdot p_2$ with $p_1=
(0,0,i)$ and $p_2=(0,0,-i)$. Consequently, every $\mbb{C}$--orbit is closed.
\end{proof}

\begin{rem} The orbit space $X_f/\mbb{C}$ is the unit disc
with a doubled origin and in particular not Hausdorff.
\end{rem}

We calculate slices at the point $p_j$, $j=1,2$, as follows. Let
$\varphi_j\colon\Delta\times\mbb{C}\to X_f$ be given by $\varphi_1(z,t):=
t\cdot(z,0,i)$ and $\varphi_2(w,s)=s\cdot(w,0,-i)$. Solving the equation
$s\cdot(w,0,-i)=t\cdot(z,0,i)$ for $(w,s)$ yields the transition function
$\varphi_{12}=\varphi_2^{-1}\circ\varphi_1\colon\Delta^*\times\mbb{C}\to
\Delta^*\times\mbb{C}$,
\begin{equation*}
(z,t)\mapsto\left(z,t+\frac{2i}{f(z)}\right).
\end{equation*}
The function $\frac{1}{f}$ is a meromorphic function on $\Delta$ without zeros
and with the unique pole $0$.

\begin{lem}
Let $\mbb{R}$ act on $X_f$ via $\mbb{R}\hookrightarrow\mbb{C}$, $t\mapsto ta$,
for some $a\in\mbb{C}^*$. Then there is no $\mbb{R}$--invariant domain
$D\subset X_f$ with $D\cap\mbb{C}\cdot p_j\not=\emptyset$ for $j=1,2$ on which
$\mbb{R}$ acts properly.
\end{lem}

\begin{proof}
Suppose that $D\subset X_f$ is an $\mbb{R}$--invariant domain with $D\cap
\mbb{C}\cdot p_j\not=\emptyset$ for $j=1,2$. Without loss of generality we may
assume that $p_1\in D$ and $\zeta\cdot p_2=(0,-2\zeta i,-i)\in D$ for some
$\zeta\in\mbb{C}$. We will show that the orbits $\mbb{R}\cdot p_1$ and $\mbb{R}
\cdot(\zeta\cdot p_2)$ cannot be separated by $\mbb{R}$--invariant open
neighborhoods.

Let $U_1\subset D$ be an $\mbb{R}$--invariant open neighborhood of $p_1$. Then
there are $r,r'>0$ such that $\Delta_r^*\times\Delta_{r'}\times\{i\}\subset
U_1$ holds. Here, $\Delta_r=\{z\in\mbb{C};\ \abs{z}<r\}$. For $(\eps_1,\eps_2)
\in\Delta_r^*\times\Delta_{r'}$ and $t\in\mbb{R}$ we have
\begin{equation*}
t\cdot(\eps_1,\eps_2,i)=\bigl(\eps_1,\eps_2+2(ta)i+(ta)^2f(\eps_1),i+(ta)
f(\eps_1)\bigr)\in U_1.
\end{equation*}
We have to show that for all $r_2,r_3>0$ there exist $(\wt{\eps_2},\wt{\eps_3})
\in\Delta_{r_2}\times\Delta_{r_3}$, $(\eps_1,\eps_2)\in\Delta^*_r\times
\Delta_{r'}$ and $t\in\mbb{R}$ such that
\begin{equation}\label{Equn:1}
\bigl(\eps_1,\eps_2+2(ta)i+(ta)^2f(\eps_1),i+(ta)f(\eps_1)\bigr)=
(\eps_1,-2\zeta i+\wt{\eps_2},-i+\wt{\eps_3})
\end{equation}
holds.

Let $r_2,r_3>0$ be given. From~\eqref{Equn:1} we obtain
$\wt{\eps_3}=taf(\eps_1)+2i$ or, equivalently, $ta=
\frac{\wt{\eps_3}-2i}{f(\eps_1)}$. Setting $\wt{\eps_2}=\eps_2$ we obtain from
$2(ta)i+(ta)^2f(\eps_1)=-2\zeta i$ the equivalent expression
\begin{equation}\label{Equn:2}
f(\eps_1)=-2i\frac{\zeta+ta}{(ta)^2}.
\end{equation}
for $t\not=0$. Choosing a real number $t\gg1$, we find an $\eps_1\in\Delta^*_r$
such that~\eqref{Equn:2} is fulfilled. After possibly enlarging $t$ we have
$\wt{\eps_3}:=taf(\eps_1)+2i=-2i\frac{\zeta}{ta}\in\Delta_{r_3}$. Together with
$\eps_2=\wt{\eps_2}$ equation~\eqref{Equn:1} is fulfilled and the proof is
finished.
\end{proof}

Thus, the Stein surface $X_f$ cannot be obtained as globalization of the local
$\mbb{C}$--action on any $\mbb{R}$--invariant domain $D\subset X_f$ on which
$\mbb{R}$ acts properly.

\subsection{The quotient $X^*/\mbb{C}$ is Hausdorff}

Suppose that $X^*/\mbb{C}$ is not Hausdorff and let $x_1,x_2\in X$ be such that
the corresponding $\mbb{C}$--orbits cannot be separated in $X^*/\mbb{C}$. Since
we already know that $\mbb{C}$ acts locally proper on $X^*$ we find local
holomorphic slices $\varphi_j\colon\Delta\times\mbb{C}\to U_j\subset X$,
$\varphi_j(z,t)=t\cdot s_j(z)$ at each $\mbb{C}\cdot x_j$ where $s_j\colon
\Delta\to X$ is holomorphic with $s_j(0)=x_j$. Consequently, we obtain the
transition function $\varphi_{12}\colon(\Delta\setminus A)\times\mbb{C}\to
(\Delta\setminus A)\times\mbb{C}$ for some closed subset $A\subset\Delta$ which
must be of the form $(z,t)\mapsto\bigl(z,t+f(z)\bigr)$ for some $f\in
\mathcal{O}(\Delta\setminus A)$. The following lemma applies to show that $A$
is discrete and that $f$ is meromorphic on $\Delta$. Hence, we are in one of
the model cases discussed in the previous subsection.

\begin{lem}
Let $\Delta_1$ and $\Delta_2$ denote two copies of the unit disk
$\{z\in\mbb{C};\ \abs{z}<1\}$. Let $U\subset\Delta_j$, $j=1,2$, be a connected
open subset and $f\colon U\subset\Delta_1\to\mbb{C}$ a non-constant holomorphic
function on $U$. Define the complex manifold
\begin{equation*}
M:=(\Delta_1\times\mbb{C})\cup(\Delta_2\times\mbb{C})/\!\!\thicksim,
\end{equation*}
where $\thicksim$ is the relation $(z_1,t_1)\thicksim(z_2,t_2):\Leftrightarrow
z_1=z_2=:z\in U$ and $t_2=t_1+f(z)$.

Suppose that $M$ is Hausdorff. Then the complement $A$ of $U$ is discrete and 
$f$ extends to a meromorphic function on $\Delta_1$.
 \end{lem}

\begin{proof}
We first prove that for every sequence $(x_n)$, $x_n\in U$, with
$\lim_{n\to\infty}x_n=p\in\partial U$, one has
$\lim_{n\to\infty}\abs{f(x_n)}=\infty\in\mbb{P}_1(\mbb{C})$. Assume the
contrary, i.e.\ there is a sequence $(x_n)$, $x_n \in U$,  with $\lim_{n\to
\infty}x_n=p\in\partial U$ such that $\lim_{n\to\infty}f(x_n)=a\in\mbb{C}$.
Choose now $t_1\in\mbb{C}$, consider the two points $(p,t_1)\in\Delta_1\times
\mbb{C}$ and $(p,t_1+a)\in\Delta_2\times\mbb{C}$ and note their corresponding
points in $M$ as $q_1$ and $q_2$. Then $q_1\not=q_2$. The sequences
$(x_n,t_1)\in\Delta_1\times\mbb{C}$ and $(x_n,t_1+f(x_n))\in\Delta_2\times
\mbb{C}$ define the same sequence in $M$ having $q_1$ and $q_2$ as accumulation
points. So $M$ is not Hausdorff, a contradiction.

In particular we have proved that the zeros of $f$ do not accumulate to 
$\partial U$ in $\Delta_1$. So there is an open neighborhood $V$ of $\partial
U$ in $\Delta_1$  such that the restriction of $f$ to $W:=U\cap V$ does not
vanish. Let $g:=1/f$  on $W$. Then $g$ extends to a continuous function on $V$
taking the value zero outside of $U$. The theorem of Rado implies that this
function is holomorphic on $V$. It follows that the boundary $\partial U$ is
discrete in $\Delta_1$ and that $f$ has a pole in each of the points of this
set, so $f$ is a meromorphic function on $\Delta_1$.
\end{proof}

\begin{thm}\label{properaction}
The orbit space $X^*/\mbb{C}$ is Hausdorff. Consequently, $\mbb{C}$ acts
properly on $X^*$.
\end{thm}

\begin{proof}
By virtue of the above lemma, in a neighborhood of two non-separable
$\mbb{C}$--orbits $X$ is isomorphic to a domain in one of the model Stein
surfaces discussed in the previous subsection. Since we have seen there that
these surfaces are never globalizations, we arrive at a contradiction. Hence,
all $\mbb{C}$--orbits are separable.
\end{proof}

\section{Examples}

In this section we discuss several examples which illustrate our results.

\subsection{Hyperbolic Stein surfaces with proper $\mbb{R}$--actions}

Let $R$ be a compact Riemann surface of genus $g\geq2$. It follows that the
universal covering of $R$ is given by the unit disc $\Delta\subset\mbb{C}$ and
hence that $R$ is hyperbolic. The fundamental group $\pi_1(R)$ of $R$ contains
a normal subgroup $N$ such that $\pi_1(R)/N\cong\mbb{Z}$. Let $\wt{R}\to R$
denote the corresponding normal covering. Then $\wt{R}$ is a hyperbolic Riemann
surface with a holomorphic $\mbb{Z}$--action such that $\wt{R}/\mbb{Z}=R$. Note
that $\mbb{Z}$ is not contained in a one parameter group of automorphisms of
$\wt{R}$.

We have two mappings
\begin{equation*}
\xymatrix{
X:=\mbb{H}\times_\mbb{Z}\wt{R}\ar[r]^{q}\ar[d]_{p} & \wt{R}/\mbb{Z}=R\\
\mbb{H}/\mbb{Z}\cong\Delta\setminus\{0\}.
}
\end{equation*}

The map $p\colon X\to\Delta\setminus\{0\}$ is a holomorphic fiber bundle with
fiber $\wt{R}$. Since the Serre problem has a positive answer if the fiber is
a non-compact Riemann surface (\cite{Mok}), the suspension
$X=\mbb{H}\times_\mbb{Z}\wt{R}$ is a hyperbolic Stein surface. The group
$\mbb{R}$ acts on $\mbb{H}\times\wt{R}$ by $t\cdot(z,x)=(z+t,x)$ and this
action commutes with the diagonal action of $\mbb{Z}$. Consequently, we obtain
an action of $\mbb{R}$ on $X$.

\begin{lem}
The universal globalization of the local $\mbb{C}$--action on $X$ is given by
$X^*=\mbb{C}\times_\mbb{Z}\wt{R}$. Moreover, $\mbb{C}$ acts properly on $X^*$.
\end{lem}

\begin{proof}
One checks directly that $t\cdot[z,x]:=[z+t,x]$ defines a holomorphic
$\mbb{C}$--action on $X^*=\mbb{C}\times_\mbb{Z}\wt{R}$ which extends the
$\mbb{R}$--action on $X$. We will show that $X$ is orbit-connected in $X^*$:
Since $[z+t,x]$ lies in $X$ if and only if there exist elements
$(z',x')\in\mbb{H}\times\wt{R}$ and $m\in\mbb{Z}$ such that $(z+t,x)=
\bigl(z'+m,m\cdot x'\bigr)$, we conclude $\mbb{C}[z,x]=\bigl\{t\in\mbb{C};\
\im(t)>-\im(z)\bigr\}$ which is connected.

In order to show that $\mbb{C}$ acts properly on $X^*$ it is sufficient to show
that $\mbb{C}\times\mbb{Z}$ acts properly on $\mbb{C}\times\wt{R}$. Hence, we
choose sequences $\{t_n\}$ in $\mbb{C}$, $\{m_n\}$ in $\mbb{Z}$ and
$\bigl\{(z_n,x_n)\bigr\}$ in $\mbb{C}\times\wt{R}$ such that
\begin{equation*}
\bigl((t_n,m_n)\cdot(z_n,x_n),(z_n,x_n)\bigr)=\bigl((z_n+t_n+m_n,m_n\cdot
x_n),(z_n,x_n)\bigr)\to\bigl((z_1,x_1),(z_0,x_0)\bigr)
\end{equation*}
holds. Since $\mbb{Z}$ acts properly on $\wt{R}$, it follows that $\{m_n\}$ has
a convergent subsequence, which in turn implies that $\{t_n\}$ has a convergent
subsequence. Hence, the lemma is proven.
\end{proof}

\begin{prop}
The quotient $X/\mbb{Z}\cong\Delta^*\times R$ is not holomorphically separable
and in particular not Stein. The quotient $X^*/\mbb{C}$ is biholomorphically
equivalent to $\wt{R}/\mbb{Z}=R$.
\end{prop}

\begin{proof}
It is sufficient to note that the map $\Phi\colon X=\mbb{H}\times_\mbb{Z}\wt{R}
\to\Delta^*\times R$, $Phi[z,x]:=\bigl(e^{2\pi iz},[x]\bigr)$, induces a
biholomorphic map $X/\mbb{Z}\to\Delta^*\times R$.
\end{proof}

\begin{prop}
The quotient $X/\mbb{Z}\cong\Delta^*\times R$ is not holomorphically separable
and in particular not Stein.
\end{prop}

Thus we have found an example for a hyperbolic Stein surface $X$ endowed with a
proper $\mbb{R}$--action such that the associated $\mbb{Z}$--quotient is not
holomorphically separable. Moreover, the $\mbb{R}$--action on $X$ extends to a
proper $\mbb{C}$--action on a Stein manifold $X^*$ containing $X$ as an
orbit-connected domain such that $X^*/\mbb{C}$ is any given compact Riemann
surface of genus $g\geq2$.

\subsection{Counterexamples with domains in $\mbb{C}^n$}

There is a bounded Reinhardt domain $D$ in $\mbb{C}^2$ endowed with a
holomorphic action of $\mbb{Z}$ such that $D/\mbb{Z}$ is not Stein. However,
this $\mbb{Z}$--action does not extend to an $\mbb{R}$--action. We give
quickly the construction.

Let $\lambda := \frac{1}{2}(3+\sqrt{5})$ and $$D:=\{(x,y) \in \mathbb C^2 \mid
\vert x\vert > \vert y \vert^{\lambda}, \vert y\vert > \vert x \vert^{\lambda}
\}.$$
It is obvious that $D$ is a bounded Reinhardt domain in $\mathbb C^2$ avoiding
the coordinate hyperplanes. The holomorphic automorphism group of $D$ is a
semidirect product $\Gamma \ltimes (S^1)^2$, where the group $\Gamma \simeq
\mathbb Z$ is generated by the automorphism $(x,y) \mapsto (x^3y^{-1},x)$ and
$(S^1)^2$ is the rotation group. Therefore the group $\Gamma$ is not contained
in a one-parameter group. Furthermore the quotient $D/\Gamma $ is the
(non-Stein) complement of the singular point in a $2$-dimensional normal
complex Stein space, a so-called "cusp singularity". These singularities are
intensively studied in connection with Hilbert modular surfaces
and Inoue-Hirzebruch surfaces, see e.g. \cite{vdG} and \cite{Zaff}.

In the rest of this subsection we give an example of a hyperbolic domain of
holomorphy in a $3$--dimensional Stein solvmanifold endowed with a proper
$\mbb{R}$--action such that the $\mbb{Z}$--quotient is not Stein. While this
domain is not simply-connected, its fundamental group is much simpler than the
fundamental groups of our two-dimensional examples.

Let $G:=\left\{\left(\begin{smallmatrix}1&a&c\\0&1&b\\0&0&1\end{smallmatrix}
\right);\ a,b,c\in\mbb{C}\right\}$ be the complex Heisenberg group and let us
consider its discrete subgroup
\begin{equation*}
\Gamma:=\left\{\begin{pmatrix}1&m&\frac{m^2}{2}+2\pi ik\\0&1&m+2\pi il\\0&0&1
\end{pmatrix};\ m,k,l\in\mbb{Z}\right\}.
\end{equation*}
Note that $\Gamma$ is isomorphic to $\mbb{Z}_m\ltimes\mbb{Z}^2_{(k,l)}$. We let $\Gamma$ act
on $\mbb{C}^2$ by
\begin{equation*}
(z,w)\mapsto\left(z+mw-\frac{m^2}{2}-2\pi ik,w-m-2\pi il\right).
\end{equation*}

\begin{prop}
The group $\Gamma$ acts properly and freely on $\mbb{C}^2$, and the quotient
manifold $\mbb{C}^2/\Gamma$ is holomorphically separable but not Stein.
\end{prop}

\begin{proof}
Since $\Gamma'\cong\mbb{Z}^2$ is a normal subgroup of $\Gamma$, we obtain
$\mbb{C}^2/\Gamma\cong(\mbb{C}^2/\Gamma')/(\Gamma/\Gamma')$. The map $\mbb{C}^2
\to\mbb{C}^*\times\mbb{C}^*$, $(z,w)\mapsto\bigl(\exp(z),\exp(w)\bigr)$,
identifies $\mbb{C}^2/\Gamma'$ with $\mbb{C}^*\times\mbb{C}^*$. The induced
action of $\Gamma/\Gamma'\cong\mbb{Z}$ on $\mbb{C}^*\times\mbb{C}^*$ is given
by
\begin{equation*}
(z,w)\mapsto\left(e^{-m^2/2}zw^m,e^{-m}w\right)
\end{equation*}
which shows that $\Gamma$ acts properly and freely on $\mbb{C}^2$. Moreover, we
obtain the commutative diagram
\begin{equation*}
\xymatrix{
\mbb{C}^*\times\mbb{C}^*\ar[d]_{(z,w)\mapsto w}\ar[r] &
Y:=(\mbb{C}^*\times\mbb{C}^*)/\mbb{Z}\ar[d]\\
\mbb{C}^*\ar[r] & T:=\mbb{C}^*/\mbb{Z}.
}
\end{equation*}

The group $\mbb{C}^*$ acts by multiplication in the first factor on
$\mbb{C}^*\times\mbb{C}^*$ and this action commutes with the $\mbb{Z}$--action.
One checks directly that the joint $(\mbb{C}^*\times\mbb{Z})$--action on
$\mbb{C}^*\times\mbb{C}^*$ is proper which implies that the map $Y\to T$ is
a $\mbb{C}^*$--principal bundle. Conseqently, $Y$ is not Stein.

In order to show that $Y$ is holomorphically separable, note that by~\cite{Oel}
this $\mbb{C}^*$--principal bundle $Y\to T$ extends to a line bundle $p\colon
L\to T$ with first Chern class $c_1(L)=-1$. Therefore the zero section of
$p\colon L\to T$ can be blown down and we obtain a singular normal Stein space
$\overline{Y}=Y\cup\{y_0\}$ where $y_0=\Sing(\overline{Y})$ is the blown down
zero section. Thus $Y$ is holomorphically separable.
\end{proof}

Let us now choose a neighborhood of the singularity $y_0\in\ol{Y}$
biholomorphic to the unit ball and let $U$ be its inverse image in $\mbb{C}^2$.
It follows that $U$ is a hyperbolic domain with smooth strictly Levi-convex
boundary in $\mbb{C}^2$ and in particular Stein. In order to obtain a proper
action of $\mbb{R}$ we form the suspension $D=\mbb{H}\times_\Gamma U$ where
$\Gamma$ acts on $\mbb{H}\times U$ by
$(t,z,w)\mapsto(t+m,z+mw-\frac{m^2}{2}-2\pi ik,w-m-2\pi il)$.

\begin{prop}
The suspension $D=\mbb{H}\times_\Gamma U$ is isomorphic to a Stein domain in
the Stein manifold $G/\Gamma$.
\end{prop}

\begin{proof}
We identify $\mbb{H}\times U$ with the $\mbb{R}\times\Gamma$--invariant domain
\begin{equation*}
\Omega:=\left\{\begin{pmatrix}1&a&c\\0&1&b\\0&0&1\end{pmatrix};\ \im(a)>0,
(c,b)\in U\right\}
\end{equation*}
in $G$.

Since $\mbb{H}\times U$ is Stein, it follows that $\mbb{H}\times_\Gamma U$ is
locally Stein in $G/\Gamma$. Hence, by virtue of~\cite{DocGr} we only have to
show that $G/\Gamma$ is Stein.

For this we note first that $G$ is a closed subgroup of
${\rm{SL}}(2,\mbb{C})\ltimes\mbb{C}^2$ which implies that $G/\Gamma$ is a
closed complex submanifold of
$X:=\bigl({\rm{SL}}(2,\mbb{C})\ltimes\mbb{C}^2\bigr)/\Gamma$. By~\cite{Oel} the
manifold $X$ is holomorphically separable, hence $G/\Gamma$ is holomorphically
separable. Since $G$ is solvable, a result of Huckleberry and Oeljeklaus
(\cite{HuOe2}) yields the Steinness of $G/\Gamma$.

One checks directly that the action of $\mbb{R}\times\Gamma$ on $\mbb{H}\times
U$ is proper which implies that $\mbb{R}$ acts properly on
$\mbb{H}\times_\Gamma U$.
\end{proof}

Because of $(\mbb{H}\times_\Gamma U)/\mbb{Z}\cong \Delta^*\times (U/\Gamma)$
this quotient manifold is not Stein but holomorphically separable.

\section{Bounded domains with proper $\mbb{R}$--actions}

In this section we give the proof of our main result.

\subsection{Proper $\mbb{R}$--actions on $D$}

Let $D\subset\mbb{C}^n$ be a bounded domain and let $\Aut(D)^0$ be the
connected component of the identity in $\Aut(D)$.

\begin{lem}
A proper $\mbb{R}$--action by holomorphic transformations on $D$ exists if and
only if the group $\Aut(D)^0$ is non-compact.
\end{lem}

The proof follows from the existence of a diffeomorphism $K\times
V\to\Aut(D)^0$ where $K$ is a maximal compact subgroup of $\Aut(D)^0$ and $V$
is a linear subspace of the Lie algebra of $\Aut(D)^0$.

\subsection{Steinness of $D/\mbb{Z}$}
Now we give the proof of our main result.

\begin{thm}
Let $D$ be a simply-connected bounded domain of holomorphy in $\mbb{C}^2$.
Suppose that the group $\mbb{R}$ acts properly by holomorphic transformations
on $D$. Then the complex manifold $D/\mbb{Z}$ is biholomorphically equivalent
to a domain of holomorphy in $\mbb{C}^2$.
\end{thm}

\begin{proof}
Let $D\subset\mbb{C}^2$ be a simply-connected bounded domain of holomorphy.
Since the Serre problem is solvable if the fiber is $D$, see \cite{Siu}, the
universal globalization $D^*$ is a simply-connected Stein surface,
\cite{IanTaTr}. Moreover, we have shown in Theorem~\ref{properaction}, that
$\mbb{C}$ acts properly on $D^*$. Since the Riemann surface
$D^*/\mbb{C}$ is also simply-connected, it must be $\Delta$, $\mbb{C}$ or
$\mbb{P}_1(\mbb{C})$. In all three cases the bundle $D^*\to D^*/\mbb{C}$
is holomorphically trivial. So we can exclude the case that $D^*/\mbb{C}$ is
compact and it follows that $D/\mbb{Z}\cong\mbb{C}^*\times(D^*/\mbb{C})$ is a
Stein domain in $\mbb{C}^2$.
\end{proof}

\subsection{A normal form for domains with non-compact $\Aut(D)^0$}

Let $D\subset\mbb{C}^2$ be a simply-connected bounded domain of holomorphy such
that the identity component of its automorphism group is non-compact. As we
have seen, this yields a proper $\mbb{R}$--action on $D$ by holomorphic
transformations and the universal globalization of the induced local
$\mbb{C}$--action on $D$ is isomorphic to $\mbb{C}\times S$ where $S$ is either
$\Delta$ or $\mbb{C}$ and where $\mbb{C}$ acts by translation in the first
factor.

Moreover, there are plurisubharmonic functions $u,-o\colon\mbb{C}\times S\to
\mbb{R}\cup\{-\infty\}$ which fulfill
\begin{equation*}
u\bigl(t\cdot(z_1,z_2)\bigr)=u(z_1,z_2)-\im(t)\quad\text{and}\quad
o\bigl(t\cdot(z_1,z_2)\bigr)=o(z_1,z_2)-\im(t)
\end{equation*}
such that $D=\bigl\{(z_1,z_2)\in\mbb{C}\times S;\ u(z_1,z_2)<0<o(z_1,z_2)
\bigr\}$. From this we conclude $u(z_1,z_2)=u(0,z_2)-\im(z_1)$,
$o(z_1,z_2)=o(0,z_2)-\im(z_1)$ and define $u'(z_2):=u(0,z_2)$,
$o'(z_2):=o(0,z_2)$.

We summarize our remarks in the following

\begin{thm}
Let $D$ be a simply-connected bounded domain of holomorphy in $\mbb{C}^2$
admitting a non-compact connected identity component of its automorphism group.
Then $D$ is biholomorphic to a domain of the form
\begin{equation*}
\wt{D}=\bigl\{(z_1,z_2)\in\mbb{C}\times S;\ u'(z_2)<\im(z_1)<o'(z_2)\bigr\},
\end{equation*}
where the functions $u',-o'$ are subharmonic in $S$.
\end{thm}

\end{document}